\documentclass[11pt]{article}

\usepackage{amsmath, amsthm}%
\usepackage{amsfonts}%
\usepackage{amssymb}%
\usepackage{graphicx}
\usepackage{amscd, amstext, hyperref}
\usepackage{pictexwd}
\usepackage{dcpic}
\usepackage[mathscr]{eucal}
\usepackage{multicol}

\newcommand{\arxiv}{\texttt{arXiv} }

\theoremstyle{definition}



\begin{document}

\title{Women's Representation in Mathematics Subfields\\Evidence from the \texttt{arXiv}}
\author{Abra Brisbin and Ursula Whitcher}

\date{}
\maketitle

%\abstract{}

\section{Women in Mathematics and the Sciences}

A few years ago, while applying to run a conference for undergraduates, one of the authors wrote, ``Mathematical physics lags behind other subfields of mathematics in participation of women.''  Our grants coordinator asked for a citation; to our surprise, we found that only very general statistics for the participation of women in research mathematics, along the lines of ``algebra'' versus ``geometry'', were available.  Nevertheless, mathematical folklore is full of more specific speculation about women's research participation, such as ``Women show a distinct preference for discrete mathematics over analytical fields'' (reported in \cite{murray}).  In this essay, we provide a detailed comparison of the fraction of women working in different subfields of mathematics, using data from the \arxiv preprints server.  Our measurements are a first effort at testing hypotheses about the participation of women in subfields of mathematics.  We take the hypothesis that subfields linked to physics have fewer women as a particular case study.

Women make up approximately half the population of the United States, and more than half of all college students: in the academic year 2011-2012, for example, 56.5\% of enrolled US undergraduates and 57.4\% of bachelor's degree recipients were female.  On the other hand, women received only 43.1\% of the bachelor's degrees awarded in mathematics and statistics in the same year. (See \cite{nsf} for more statistics on degrees and enrollment.)  At higher levels, the discrepancies are starker. Recent surveys by the American Mathematical Society show that in 2012-2013, 31.2\% of the people receiving Ph.D.s were women, and women held 26.6\% of all faculty jobs in mathematics and statistics (see \cite{amsphds,amsjobs}).  Because the fraction of women studying and working on mathematics is not representative of the population as a whole, education and labor force researchers say women are \emph{underrepresented} in mathematics.

Underrepresentation of women is hardly unique to math: it is a common phenomenon across the STEM (science, technology, engineering, and mathematics) fields, as well as philosophy and economics.  We give the fraction of women receiving Ph.D.s in different branches of science in Table~\ref{Ta:sciencephds}, using recent National Science Foundation data, for the academic year 2011-2012 (\cite{nsf}). The fraction of women participating in academic fields varies greatly.  From Table~\ref{Ta:sciencephds}, we see that the fraction of female Ph.D. recipients in chemistry is nearly twice as high as the fraction of female Ph.D. recipients in physics, for example.

\begin{table}[!htbp]
\centering

\begin{tabular}{|c|c|}
\hline
\textbf{Scientific Field} & \textbf{\% Women, 2012} \\
\hline
Biological Sciences & 53.1\% \\
Earth, Atmospheric, and Ocean Sciences & 43.3\%\\
Chemistry & 39.1\%\\
Mathematics and Statistics & 28.2\%\\
Computer Sciences & 21.4\%\\
Physics & 20.0\%\\
\hline
\end{tabular}

\caption{Doctoral Degrees Awarded to Women, 2011-2012}\label{Ta:sciencephds}
\end{table}

Fields with greater representation of women are likely to have support structures which benefit women, such as child care at conferences, because more women advocate for those support structures.  The National Academy of Sciences' Committee on Maximizing the Potential of Women in Academic Science and Engineering (\cite{beyondbias}, p. 187) describe a ``social tipping point" that occurs when women's participation in a field reaches approximately 20\%, at which point women begin appearing in leadership positions and working together for support structures which benefit them.  These support structures also benefit men: for example, men with children take advantage of conference child care.  In addition to their immediate practical benefits, these support structures may allow women in these fields to devote more mental energy to their work, allowing them to produce more and better work than women in fields where they are less represented.  Women in fields with lower representation of women may experience greater stereotype threat.  Stereotype threat refers to the awareness that poor performance may reinforce others' negative stereotypes about your group.  It creates added pressure to perform, and tends to decrease performance (\cite{ssq}).

Research shows that diverse teams are more innovative. In academia, Freeman and Huang (\cite{fh}) have shown that papers with ethnically diverse co-authors are published in higher-impact journals and receive more citations (even controlling for past publishing records of the authors).  In industry, male groups of engineers have designed products that fail, sometimes catastrophically, when women use them.  Examples range from voice-recognition systems that could not process female voices, to airbags designed for the average male body that injured women when inflating (see \cite{mf}).  Thus, recruiting more women to STEM fields may be advantageous for the fields as a whole.

The authors of \cite{LCMF} showed that underrepresentation of women among Ph.D. recipients in a field is highly correlated with the belief that brilliance is necessary to succeed in that field.  Furthermore, belief that a field requires brilliance was a better predictor of the fraction of women in that field than other plausible hypotheses.  For example, though one might hypothesize that women prefer fields of research with flexible hours to fields that require long hours working in a scientific laboratory, this would not explain the low numbers of women in computer science and comparatively high numbers of women in chemistry.  In the \cite{LCMF} study, mathematicians were more likely than academics in other STEM fields to agree with statements such as ``Being a top scholar of [discipline] requires a special aptitude that just can't be taught''.  Indeed, the only academic field whose members thought it demanded more brilliance than mathematics was philosophy.

Nosek and Smyth \cite{SN15} used the Implicit Attitudes Test, a test of response speed in categorizing terms, to investigate people's implicit or unconscious associations between science and gender.  They found that women in STEM fields demonstrate lower levels of implicit science-is-male stereotyping than women in non-STEM fields, and suggest that a strong association between oneself and science could partially counteract the effects of cultural stereotypes.  Another possibility is that lower levels of implicit stereotyping may contribute to women's interest, or expectations of success, in STEM fields.  Nosek and Smyth also found lower levels of implicit stereotyping among women with PhDs in the physical sciences and engineering, compared to women with Masters degrees in the same fields.  However, the opposite association held among women in the biological sciences, where women are less underrepresented.  This could indicate that women with higher levels of implicit stereotyping are less likely to pursue doctoral degrees in fields where women are underrepresented.

Not only does the representation of women vary between STEM fields, it varies within the field of mathematics and statistics.  Women received 44.5\% of the Ph.D.s in statistics and biostatistics in 2012-2013, for example, but only 17.0\% of the Ph.D.s in analysis.  We wish to develop methods to measure the participation of women in different subfields of mathematics.  Subfields that are particularly successful in recruiting, retaining, and fostering research productivity of women might provide recruitment models for other fields of mathematics; conversely, women in subfields with very low participation of women might benefit from extra support.  Furthermore, many existing efforts to promote the participation of women in mathematics are subfield-specific.  For example, the Association for Women in Mathematics has received an NSF ADVANCE grant to sponsor Research Collaboration Conferences for Women in specific subfields such as number theory and algebraic combinatorics (see \cite{awmadvance}).  Measurements of the research productivity of women in specific subfields are important for directing these efforts and assessing their impact.

We may use the American Mathematics Society's annual survey of new Ph.D.s (\cite{amsphds}) to gain a general picture of the participation of women in different fields of mathematics.  The AMS data for the academic year 2012-2013 is summarized in Table~\ref{Ta:amsphds}.  We see immediately that women are comparatively well represented in applied subfields, including applied mathematics, statistics and biostatistics, and optimization, and badly represented in analysis and probability.  If we compute the representation of women in pure mathematics using the AMS Ph.D. data, omitting statistics and biostatistics, applied math, and math education, we see that women represent approximately $23.4\%$ of new pure mathematics Ph.D.s.  Our computations show that the fraction of women receiving Ph.D.s in pure mathematics is comparable to the fraction of women receiving Ph.D.s in physics and computer science, and is significantly lower than the fraction of women receiving Ph.D.s in the biological sciences.

The greater participation of women in statistics and biostatistics may be driven by the large number of biological applications of statistics.  The \cite{LCMF} survey showed that statisticians place a lower emphasis on brilliance than mathematicians; thus, statistics may seem more accessible to women.  Similarly, the high fraction of women in math education may be related to the large number of women interested in education more generally, or to the perception that education requires less brilliance than pure mathematics.

\begin{table}[!htbp]
\centering
\begin{tabular}{| c | c | c | c |}
\hline
\textbf{Field of Mathematics} & \textbf{Women} & \textbf{Men} & \textbf{\% Women} \\
\hline
Algebra/ Number Theory & 61 & 197 & 23.6\% \\
Analysis & 15 & 73 & 17.0\% \\
Geometry/ Topology & 40 & 134 & 23.0\% \\
Combinatorics/ Logic & 36 & 102 & 26.1\% \\
Probability  & 15 & 69 & 17.9\% \\
Statistics  \& Biostatistics & 255 & 318 & 44.5\% \\
Applied Math  & 66 & 148 & 30.8\% \\
Numerical Analysis  & 22 & 72 & 23.4\% \\
Optimization  & 11 & 13 & 45.8\% \\
Differential Equations  & 37  & 105 & 26.1\% \\
Math Education  & 14 &  9 & 60.9\% \\
Other/ Unknown  & 5 & 26 & 16.1\% \\
\hline
Total  & 577 & 1266 & 31.3\% \\
\hline

\end{tabular}
\caption{AMS Survey of Ph.D.s, 2012-2013}\label{Ta:amsphds}
\end{table}

Though intriguing, the AMS data is not fine-grained enough to detect differences in the participation of women in pure mathematics.  For example, commutative algebraists and analytic number theorists will take different courses in graduate school, attend different conferences, and publish in different journals, but recent Ph.D.s from both groups are classified as Algebra/ Number Theory in the AMS statistics.  We use data from the \texttt{arXiv} preprints server (\url{http://arXiv.org/}) to compare the participation of men and women in different mathematical subfields.  This measures a somewhat different metric of participation than the AMS data:  The AMS data focuses on the fields to which new Ph.D.s ``belong", on the basis of their research focus in graduate school, while data from the \texttt{arXiv} combines information about representation or ``belonging" in a subfield with information about research productivity, as measured by number of manuscripts posted to the \texttt{arXiv}.

\section{Mathematics and the \arxiv}\label{S:arxivintro}

The \texttt{arXiv} hosts electronic preprints, or ``e-prints'', in physics, mathematics, computer science, quantitative biology, quantitative finance, and statistics.  The roots of the \texttt{arXiv} lie in a 1989 string theory conference; after the conference, the astrophysicist Joanne Cohn emailed related papers to a group of interested scientists.  Her mailing list quickly grew.  In 1991, Paul Ginsparg created the website and interface for the \texttt{arXiv}, which he hosted at his workplace, the Los Alamos National Labs.  Since 2001, Ginsparg and the \texttt{arXiv} have been based at Cornell University.  (For more comments on the history of the \texttt{arXiv}, see \cite{arXivhistory}).

In its early days, the \texttt{arXiv} focused on theoretical physics papers.  It became an official repository for mathematics papers in 1995.  To this day, physicists and mathematicians are more likely to post papers to the \texttt{arXiv} than scientists in other disciplines \cite{eprints}.  In a survey of 584 authors with papers indexed by the Web of Science's Science Citation Index in the subject areas ``Mathematics", ``Mathematics, Applied", and ``Statistics \& Probability", Kristine Fowler found that $56\%$ of mathematicians have posted at least one paper to the \texttt{arXiv}, with about $30\%$ of mathematicians routinely posting preprints \cite{fowler}.  These results may give extra weight to tenured faculty, who were more likely to retain the same email address between when their paper was published and when they were contacted for the survey.  Many mathematicians who post preprints to the \texttt{arXiv} cite the early dissemination of research findings, and the better availability and visibility of both published and unpublished work, as motivations for using the \texttt{arXiv}.  Others don't see a reason to post to the \texttt{arXiv}, due to their satisfaction with traditional publication methods or due to uncertainty about whether posting to the \texttt{arXiv} is allowed by the journals in which their papers are published \cite{fowler}.

One potential concern for authors about posting preprints on the \arxiv is the lack of peer review.  This may be a particular concern in the General Mathematics category, which contains a disproportionate number of claims by amateur mathematicians.  However, authors who make use of the \arxiv are typically not worried about the lack of peer review, as errors could be detected by any of hundreds of people who receive a paper's abstract through the \texttt{arXiv}'s email notification service \cite{Jack02}.

We use \arxiv postings to analyze mathematicians' participation in different subfields of mathematics.  Counting postings to the \arxiv measures how frequently mathematicians create and share research.  Of course, research is only one of a mathematician's duties: teaching, institutional service, mentoring, outreach, and service to the profession are all important parts of mathematicians' professional identities.  Depending on an individual's interests or the requirements of a specific job, teaching or service may demand more time and attention than research.  We assume that active mathematics researchers share papers based on their research from time to time; for this reason, we believe \arxiv postings are a good proxy for the actual membership of a subfield.  However, inactive researchers might still identify with a subfield, and could have impact on its development by teaching related courses and mentoring emerging researchers; because our measurements are based on sharing of papers, they may not detect these more personal affiliations.

Of course, productivity as measured by rate of posting or publishing papers is not the only way to assess the impact of a researcher's work.  Other measurements include citation rates, grant funding, journal editorships, and invited talks.  Analyses of such measures have detected differences in the way men and women in mathematics are recognized.  For example, Greg Martin argues in \cite{martin} that low rates of women speakers at the 2014 International Congress of Mathematicians and Joint Mathematics Meetings reflects bias in the training and evaluation of mathematicians; similarly, Topaz and Sen suggest in \cite{topazsen} that the low representation of women on editorial boards of mathematics journals both reflects and maintains inequities in the profession.

In our analysis, we use a dataset of \texttt{arXiv} papers collected by the data scientist Emma Pierson.  Pierson collected 938,301 papers posted to the \texttt{arXiv} between its creation and July 2014.  She used a list of more than 40,000 names classified by native speakers to infer the gender of each author, based on given name \cite{pierson}.   Papers on the \arxiv may be cross-posted to several categories.  We extracted the papers where the primary or first category included Mathematics, such as Mathematics - Geometric Topology or Mathematics - Dynamical Systems, and classified these papers according to the subfield of their primary category.  This resulted in 174,074 mathematics papers.

The number of papers in each mathematics category in our dataset is listed in Table~\ref{Ta:subfields}.  Mathematicians mention standard practice in their fields as a reason both for and against posting papers to the \arxiv \cite{fowler}.  This is a likely reason for differing numbers of papers from different fields.  We anticipate that this factor will influence men and women in a given field equally, so standard practice in a given field will not affect the proportion of papers on the \arxiv with female authors.  The most popular category is Algebraic Geometry; the prevalence of algebraic geometry papers on the \arxiv may be increased due to the popularity of the \arxiv among physicists, and the many connections between algebraic geometry and physics research.  The \arxiv mathematics categories also include some entries such as Quantum Algebra with clear theoretical physics origins.

\begin{table}[!htbp]
\centering
\begin{tabular}{|c|c|c|c|}
\hline
\textbf{Subfield} & \textbf{\# of Papers} & \textbf{Subfield} & \textbf{\# of Papers} \\
\hline
Algebraic Geometry    &         16858 &Numerical Analysis       &       4388 \\
Probability           &         14180 & Complex Variables       &        4075\\
Combinatorics          &        13546 & Rings and Algebras      &        4069\\
Differential Geometry  &        13141 & Operator Algebras       &        3909\\
Analysis of PDEs       &        12437 & Algebraic Topology     &         3461\\
Number Theory         &         10971 & Commutative Algebra    &         3336\\
Dynamical Systems    &           7532 & Logic                  &         3158\\
Functional Analysis  &           7079 & Symplectic Geometry   &          2253\\
Geometric Topology   &           7051 & Metric Geometry       &          2154\\
Quantum Algebra      &           6651 & Spectral Theory       &          1910\\
Representation Theory  &         5896 & General Mathematics   &          1438\\
Classical Analysis, ODEs  &   5486 & Category Theory       &          1263\\
Group Theory              &      5357 & K-Theory \& Homology &          1225\\
Statistics Theory        &       4786 & General Topology      &          1131\\
Optimization and Control  &      4425 & History and Overview   &          908\\

\hline
\end{tabular}
\caption{Math Papers on the \arxiv by Primary Category}\label{Ta:subfields}
\end{table}

The gender of many of the authors in our dataset could not be identified.  This happens when a given name is used by both men and women, when a given name is not in the list of 40,000 names, or when an author uses only his or her initials.  We summarize our classification of authors by gender in Table~\ref{Ta:knowngenders}.  In all, authors identifiable as women make up $10.1\%$ of all authors in our dataset, or $17.0\%$ of authors whose gender is known.  We also report the gender of \emph{paper-authors}, sometimes known as \emph{authorships}.  For this measure, we count each author on a paper separately (for example, this paper has two paper-authors).

\begin{table}[!hbtp]
\centering
\begin{tabular}{|c|c|c|}
\hline
\textbf{Gender} & \textbf{Authors} & \textbf{Paper-Authors}\\
\hline
Men &    49337 & 178839\\
Unknown &    40746 & 110828\\
Women &    10100 & 25882\\
\hline
\end{tabular}
\caption{Gender Classification}\label{Ta:knowngenders}
\end{table}

\begin{figure}[!htbp]
\centering
\scalebox{1}{\includegraphics{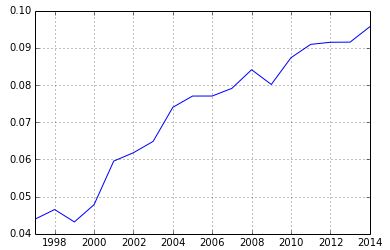}}
\caption{Women as a fraction of \texttt{arXiv} math paper-authors}\label{F:womenovertime}
\end{figure}

The participation of women in the mathematics section of the \arxiv has been increasing over time.  We illustrate the growth in the number of female mathematics paper-authors as a fraction of all mathematics paper-authors over time in Figure~\ref{F:womenovertime}.  We use the date a paper was first posted to the \arxiv.  We have omitted papers backdated to before 1995, when the \arxiv officially began accepting mathematics submissions.

According to a 2013 AMS survey, women hold $29\%$ of the full-time positions in mathematics departments.  However, the distribution is uneven: women hold only $22\%$ of the full-time positions in mathematics departments granting a doctoral degree, but $35\%$ of the full-time positions in math departments where the highest degree is a bachelor's or a master's \cite{amsjobs}.  Thus, the participation of women on the \arxiv is closer to the representation of women in doctoral mathematics departments than to the membership of the field as a whole.  Because women hold fewer positions in research-intensive departments, women may have fewer incentives than men to publish large numbers of research papers.  They may also have more difficulty accessing resources, such as funding for conference travel or course releases, that support the production and publication of research.  Alternatively, women who are highly productive early in their careers might be simultaneously more likely to post to the \arxiv and to attain full-time positions in doctoral departments.

Women might also choose to post papers using only their initials, in order to avoid gender bias.  If so, our measure of women on the \arxiv would underestimate their true participation.  Emma Pierson's analysis of \arxiv papers across scientific fields shows that women may be more likely than men to use only initials when posting to the \arxiv (see \cite{piersoninitials}).  The use of initials rather than full names also varies by subfield; we will explore this phenomenon in more detail later in this essay.

We may compare our study of gender and mathematics on the \arxiv to an analysis of gender and publishing patterns using publications recorded in the zbMATH database, which appeared while this manuscript was undergoing editorial review \cite{zbmath}.  The authors of \cite{zbmath} studied 2,245,205 papers published in ``core math journals'' since 1970.  They focused on authors of mathematics papers who could be uniquely identified, excluding, for example, combinations of initials and last names that could plausibly have been used by multiple people.  They used a list of approximately 42,000 given names to assign gender to those uniquely identified authors whose given names were known.  The fraction of ambiguous names was comparable to the fraction in our study: they write, ``We were able to assign a gender to ∼55\% of all profiles with a real first name; among those, 27,596 (∼19\%) authors were classified as women and 116,657 as men.''  The authors of \cite{zbmath} found a steady increase in the fraction of women paper-authors since 1970, which is consistent with our measurement of increasing women \arxiv paper-authors since 1995.

\section{Ranking Subfields by Paper-Authors}

The average paper in our dataset has 0.149 female authors, 1.027 male authors, and 0.637 authors
whose gender could not be determined automatically, for a total of 1.813 total authors.  We rank \arxiv categories by average number of female paper-authors in Table~\ref{Ta:paperauthor}.

\begin{table}[!htbp]
\centering
\begin{tabular}{|c|c|c|c|c|}
\hline
\textbf{Field} & \textbf{Women} & \textbf{Men} & \textbf{Unknown} & \textbf{Total} \\
\hline

General Mathematics	&	0.041 &	0.658	& 0.682 &	1.381 \\
Quantum Algebra	&	0.086 &	0.831&	0.857&	1.774\\
Operator Algebras	&	0.092&	1.102	&0.532&	1.73\\
K-Theory and Homology	&	0.096&	1.049&	0.496&	1.642\\
Spectral Theory	&	0.102	&1.163	&0.594	&1.860\\
Complex Variables	&	0.103&	0.975&	0.653&	1.731\\
Number Theory	&	0.105&	0.950&	0.544&	1.599\\
Logic	&	0.107	&1.001&	0.519&	1.627\\
Category Theory	&	0.117&	0.952&	0.463&	1.532\\
Differential Geometry	&	0.118	&0.938	&0.650	&1.707\\
History and Overview	&	0.120 &	0.939&	0.304	&1.363\\
General Topology	&	0.132&	0.916&	0.609	&1.657 \\
Metric Geometry	&	0.134 &	1.105&	0.508&	1.747 \\
Algebraic Geometry	&	0.140 &	0.950 &	0.535 &	1.624 \\
Classical Analysis and ODEs	&	0.141&	0.968&	0.702&	1.812\\
Algebraic Topology	&	0.142&	1.053&	0.442&	1.637\\
Symplectic Geometry	&0.143	&0.907&	0.590 &	1.640\\
Functional Analysis	&	0.145 &	1.019 &	0.697 &	1.861 \\
Dynamical Systems	&	0.150 &	1.079 &	0.652 &	1.882\\
Probability	&	0.150 &	1.165 &	0.649&	1.964 \\
Representation Theory	&	0.153&	0.975	&0.534&	1.661 \\
Geometric Topology	&	0.167 &	0.951&	0.532 &	1.651\\
Group Theory	&0.169	&1.074	&0.517 &	1.760 \\
Optimization and Control	&	0.173 &	1.289 &	0.828 &	2.289 \\
Rings and Algebras	&	0.174	&0.833&	0.765	&1.772\\
Numerical Analysis &	0.188 	&1.225 	&0.907 	&2.321 \\
Analysis of PDEs	&	0.193 &	1.076 &	0.745 &	2.014 \\
Combinatorics	&	0.202 &	1.143 &	0.647 &	1.993 \\
Statistics Theory	&	0.205	&1.010 &	0.792 &	2.097 \\
Commutative Algebra	&	0.239 &	1.048 &	0.552 &	1.839 \\

\hline

\end{tabular}
\caption{Average Number of Paper-Authors}\label{Ta:paperauthor}
\end{table}

We see immediately that women are especially unlikely to post in the General Mathematics category.  General Mathematics contains an unusually high number of amateur claims to have proved famous conjectures; recent submissions include ``A New Way to Proof 3x+1 Problem'' and ``The topological proof of the Poincare conjecture'' (see \cite{collatz,poincare}).  Thus, the underrepresentation of women in General Mathematics may not give us much information about the careers of women who are professional mathematicians.

Many of the categories in Table~\ref{Ta:paperauthor} with high average numbers of female paper-authors also have high average numbers of authors overall.  In order to correct for this phenomenon, we compute a \emph{discrepancy} score for each category, using the fact that $8.2\%$ of the total paper-authors are identifiably female:

\[\text{discrepancy} = \frac{\text{female paper-authors} - \text{expected female paper-authors}}{\text{total paper-authors}}\]

\noindent A discrepancy of 0 would indicate that a category was completely average for our dataset; categories with negative discrepancies have worse representation of women paper-authors than average, and categories with positive discrepancies have better representation than average.  We rank \arxiv mathematics categories by discrepancy in Table~\ref{Ta:discrepancy}.  For each category, we also report a 95\% confidence interval.  We computed the confidence interval using a bootstrap method with 10,000 resamplings.  This enables us to obtain confidence intervals for our novel discrepancy score, independent of assumptions about the distribution of this statistic.

\begin{table}[!htbp]

\centering
\begin{tabular}{|c|c|c|}
\hline
\textbf{Subject} & \textbf{Discrepancy} & \textbf{Confidence Interval} \\
\hline
General Mathematics	&	-0.052 & (-0.060,-0.044)\\
Quantum Algebra	&	-0.033 & (-0.037, -0.029) \\
Operator Algebras	&	-0.029 & (-0.035, -0.023)\\
Spectral Theory	&	-0.027 & (-0.035, -0.018)\\
K-Theory and Homology &	-0.023 & (-0.034, -0.012)\\
Complex Variables	&	-0.022 & (-0.029, -0.016)\\
Number Theory	&	-0.017 & (-0.021, -0.012)\\
Logic	&	-0.016 & (-0.024, -0.009)\\
Differential Geometry	&	-0.013 & (-0.017, -0.009)\\
Optimization and Control	&	-0.006 & (-0.012, -0.001)\\
Category Theory	&	-0.006 & (-0.018, 0.008)\\
Metric Geometry	&	-0.006 & (-0.015, 0.004)\\
Probability	& -0.005& (-0.009,	-0.002) \\
Functional Analysis	&	-0.004 & (-0.009, 0.001)\\
Classical Analysis and ODEs	&	-0.004 & (-0.010, 0.002)\\
General Topology &	-0.003 & (-0.015, 0.010)\\
Dynamical Systems	&	-0.002 & (-0.007, 0.003)\\
Numerical Analysis	&	-0.001 & (-0.007, 0.005)\\
Algebraic Geometry	&	0.004 & (0.000, 0.008)\\
Algebraic Topology	&	0.005 & (-0.003, 0.013)\\
Symplectic Geometry	&	0.005 & (-0.005, 0.016) \\
History and Overview	&	0.006 & (-0.011,0.023)\\
Representation Theory	&	0.010 & (0.004, 0.016)\\
Analysis of PDEs	&	0.014 & (0.010, 0.018)\\
Group Theory	&	0.014 & (0.007, 0.020)\\
Statistics Theory	& 0.016 & (0.009, 0.022) \\
Rings and Algebras	&	0.016 & (0.008, 0.024)\\
Geometric Topology	&	0.019 & (0.013, 0.025)\\
Combinatorics	&	0.019 & (0.015,	0.023)\\
Commutative Algebra	&	0.048 & (0.039, 0.057)\\
\hline

\end{tabular}
\caption{\arxiv Categories By Discrepancy}\label{Ta:discrepancy}
\end{table}

We may use our confidence intervals to divide \arxiv categories into three groups.  If the 95\% confidence interval contains 0, representation of female paper-authors in that category is average, compared to other mathematics categories on the \texttt{arxiv}; if all elements of the confidence interval are less than 0, representation is below average, and if all elements of the confidence interval are greater than 0, representation is above average.  We summarize our grouping in Table~\ref{Ta:discrepancyconfidence}.

\begin{table}[!htbp]
\centering
\begin{tabular}{|c|p{10cm}|}
\hline
\textbf{Female Paper-Authors} & \textbf{\arxiv Categories} \\
\hline

Significantly Fewer &
Complex Variables, Differential Geometry, General Mathematics, K-Theory and Homology, Logic, Operator Algebras, Optimization and Control, Probability, Quantum Algebra, Number Theory, Spectral Theory  \\

\hline

Average & Algebraic Topology, Category Theory, Classical Analysis \& ODEs, Dynamical Systems, Functional Analysis, General Topology, History and Overview, Metric Geometry, Numerical Analysis, Symplectic Geometry \\

\hline

Significantly More & Algebraic Geometry, Analysis of PDEs, Combinatorics,Commutative Algebra, Geometric Topology, Group Theory,Representation Theory, Rings and Algebras, Statistics Theory \\

\hline
\end{tabular}
\caption{Female Paper-Authors vs. Expectation}\label{Ta:discrepancyconfidence}
\end{table}

We note that the \arxiv mathematics categories with the strongest links to physics, such as Operator Algebras and Quantum Algebra, are all in the Significantly Fewer group, providing statistical support for the anecdotal observation about mathematical physics which inspired this research.  %The underrepresentation of female paper-authors in these categories may be driven by the underrepresentation of women in physics.

Because of the increasing proportion of female paper-authors over time (Figure~\ref{F:womenovertime}), we also examined the discrepancy score for each field in two separate decades, 1995-2004 and 2005-2014.  For this analysis, the expected number of female paper-authors was computed separately based on each decade's proportion of female paper-authors.  The results of the by-decade analysis largely agreed with the overall analysis.  Three fields with significantly fewer female paper-authors overall (Differential Geometry, K-Theory and Homology, and Probability) had significantly fewer female paper-authors in 2005-2014 but not in 1995-2004.  Two fields with significantly more female paper-authors overall (Rings and Algebras, and Statistics Theory) had significantly more female paper-authors in 2005-2014 but not in 1995-2004.  It should be noted that the number of mathematics papers from 1995-2004 posted on the \arxiv is much lower than the number from 2005-2014 (32,230 versus 141,844), so small sample sizes could contribute to a lack of significance.  This is particularly true for Statistics Theory and for K-Theory and Homology, which had only 147 and 216 papers, respectively, from 1995-2004.  Of the fields that were not significantly different from average in the overall analysis, only two fields were significantly different in the by-decade analysis:  Functional Analysis had significantly fewer female paper-authors in 2005-2014, and Symplectic Geometry had significantly more female paper-authors in 1995-2004.

The authors of \cite{zbmath} report publications by women in different subfields using the 2010 Mathematics Subject Classification (MSC), a categorization system maintained by the editors of \emph{Mathematical Reviews} (MathSciNet) and zbMATH.  The MSC currently contains 63 two-digit classes, corresponding to areas of pure and applied mathematics and statistics, as well as to mathematically-oriented topics in fields such as economics and computer science.  Even within pure mathematics, the \arxiv and MSC categorization schemes do not always correspond directly: for example, Symplectic Geometry is a top-level category on the \arxiv but a subset of 53 (Differential Geometry) in the MSC, while papers categorized as 22 (Lie Groups) in the MSC might be posted to either the Group Theory or Differential Geometry sections of the \texttt{arxiv}.  In general terms, the MSC is both a finer and a more conservative classification than the \arxiv system: the MSC provides more categories, and maintains separate categories for classical mathematics subjects such as functions of a single complex variable.

The authors of \cite{zbmath} report the difference between the average percentage of women paper-authors publishing in each two-digit MSC category and the average percentage in their overall dataset.  We provide this ranking in Table~\ref{Ta:mscpaperauthors}.  Because the zbMATH dataset is proprietary, we are not able to generate confidence intervals; instead, we simply divide the categories into lowest, middle and highest thirds by percentage of female paper-authors.  Despite the differences in classification schemes, we are able to detect consensus on the role of certain subfields.  For example, K-theory and complex variables have few publications by women in both our analysis and the zbMATH data, while combinatorics and commutative algebra evidence a higher rate of women paper-authors in both analyses.

\begin{table}[!htbp]
\centering
\begin{tabular}{|c|p{10cm}|}
\hline
\textbf{Female Paper-Authors} & \textbf{MSC Categories} \\
\hline
Lowest third & 83 Relativity and gravitional theory, 19	K-Theory, 81 Quantum theory, 93	Systems theory and control, 55 Algebraic topology, 74 Mechanics of deformable solids, 11	Number theory, 70 Mechanics of particles and systems, 00 General, 76 Fluid mechanics, 32 Several complex variables, 57	Manifolds and cell complexes, 40 Sequences and series, 58 Global analysis and analysis on manifolds, 22 Topological and Lie groups, 82	Statistical mechanics, 91 Game theory and social sciences, 03 Mathematical logic and foundations, 43	Abstract harmonic analysis, 90 Operations research, 60	Probability theory

\\
\hline
Middle third & 68 Computer science, 94 Information and communication, 01 History and biography, 86	Geophysics, 46	Functional analysis, 31	Potential theory, 52	Convex and discrete geometry, 78	Optics and electromagnetic theory, 47	Operator theory, 85	Astronomy and astrophysics, 14	Algebraic geometry, 65	Numerical analysis, 18	Category theory and homological algebra, 15	Linear and multilinear algebra, 37	Dynamical systems and ergodic theory, 41	Approximations and expansions,
12	Field theory, 17	Nonassociative rings and algebras, 80	Classical thermodynamics,
28	Measure and integration,
92	Biology and other natural sciences

\\
\hline
Highest third & 20	Group theory, 33	Special functions, 53	Differential geometry, 45	Integral equations,
30	Functions of a complex variable, 34	Ordinary differential equations, 42	Harmonic analysis on Euclidean spaces, 54	General topology, 97	Mathematics education,
44	Integral transforms,
16	Associative rings and algebras,
49	Calculus of variations and optimization,
26	Real functions,
35	Partial differential equations,
51	Geometry,
39	Difference and functional equations,
13	Commutative algebra,
05	Combinatorics,
06	Order and lattices,
62	Statistics,
08	General algebraic systems

\\
\hline
\end{tabular}
\caption{MSC Paper-Authors}\label{Ta:mscpaperauthors}
\end{table}

Let us now investigate the phenomenon of \arxiv authors whose gender is unknown in more detail.  We list the top twenty given names whose gender could not be determined automatically in Table~\ref{Ta:unknownnames}.  Nineteen of these names are single initials; the remaining name belongs to Saharon Shelah, an Israeli mathematician who has published more than 1000 papers. \cite{saharon}  Many of Shelah's \arxiv contributions are in the Logic category, but he has also posted papers in sections such as Group Theory, Combinatorics, and General Topology.

\begin{table}[!htbp]
\centering
\begin{tabular}{|c|c|c|c|}
\hline
\textbf{Name} & \textbf{Frequency by Paper-Author} & \textbf{Name} & \textbf{Frequency by Paper-Author}\\
\hline

A. & 2729 & V. & 889 \\
M. & 2272 & C. & 861\\
S. & 1899 & F. & 842\\

J. & 1385 & B. & 802\\
D. & 1269 & T. & 762\\
E. & 1181 & K. & 677\\

R. & 1170 & Saharon & 648\\
G. & 1047 & N. & 596\\
P. & 936 & H. & 542\\

L. & 932 & I. & 518\\

\hline
\end{tabular}
\caption{Most Popular Names of Unknown Gender}\label{Ta:unknownnames}
\end{table}

We rank \arxiv categories by the fraction of paper-authors using only their initials in Table~\ref{Ta:initials}.  We see there is a large variation, from Logic, where only about 5\% of paper-authors use their initials, to Quantum Algebra, where nearly a third of paper-authors use only their initials.

\begin{table}[!htbp]
\centering
\begin{tabular}{|c|c|c|c|}
\hline
\textbf{Subject} & \textbf{Fraction of Initials} & \textbf{Subject} & \textbf{Fraction of Initials} \\
\hline
Logic                    &      0.054 & Differential Geometry    &      0.124\\
Statistics Theory        &      0.065 & Group Theory             &      0.125\\
Geometric Topology       &      0.076 & Dynamical Systems         &     0.125\\

Symplectic Geometry     &       0.078 & Numerical Analysis         &    0.127\\
Combinatorics           &       0.078 & Algebraic Geometry          &   0.130\\
Commutative Algebra     &       0.078 & Category Theory              &  0.132\\

Representation Theory   &       0.083 & K-Theory and Homology       &   0.133\\
Algebraic Topology      &       0.084 & General Topology            &   0.160\\
History and Overview     &      0.086 & Complex Variables            &  0.165\\

Optimization and Control &      0.086 & Classical Analysis and ODEs  &  0.165\\
Number Theory       &           0.087 & General Mathematics          &  0.167\\
Analysis of PDEs    &           0.096 & Functional Analysis          &  0.173\\

Probability        &            0.098 & Rings and Algebras           & 0.187\\
Metric Geometry     &           0.106 & Spectral Theory            &   0.202\\
Operator Algebras        &      0.123 & Quantum Algebra             &   0.290\\

\hline
\end{tabular}
\caption{Fraction of Paper-Authors Using Initials}\label{Ta:initials}
\end{table}

\begin{figure}[!htbp]
\centering
\scalebox{.75}{\includegraphics{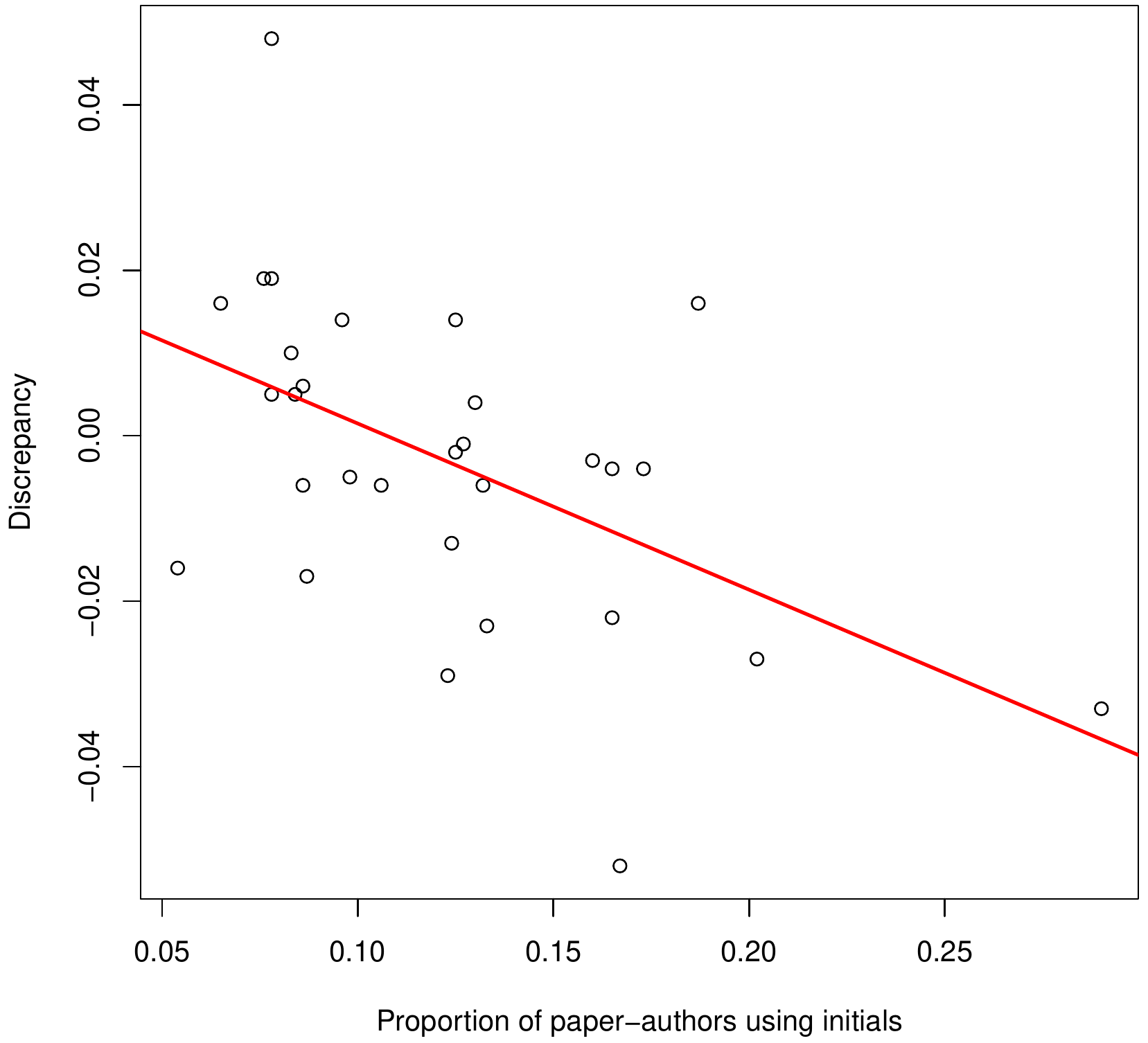}}
\caption{Initial use vs. discrepancy}\label{F:linreg}
\end{figure}

The proportion of paper-authors using only initials is associated with the discrepancy in the proportion of female paper-authors.  The correlation was -0.515 (see Figure~\ref{F:linreg}), indicating that a greater proportion of paper-authors using initials is associated with lower representation of women in the field.  If women are more likely than men to use their initials on a paper, this could partially account for the below-average representation of women observed in some subfields, such as Spectral Theory and Quantum Algebra.  However, it is also possible that below-average representation in a subfield could result in greater incentive for female authors to avoid gender bias through the use of initials.

\clearpage

\section{Ranking Subfields by Gender of Authors}

\begin{figure}[!htbp]
\centering
\scalebox{1}{\includegraphics{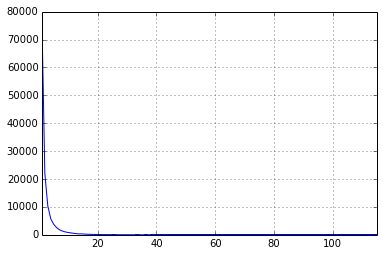}}
\caption{Paper count vs. number of authors}\label{F:papernumbers}
\end{figure}

What happens if we rank fields by the number of female authors, instead of the number of female paper-authors?  Within each \arxiv category, we identified authors with the same name.  (Note that our method fails to recognize authors who have posted papers using different names, or different combinations of given name and initials.)  The median number of papers per author within a category was 1, while the mean was 2.41 papers; the greatest number of papers per author in any single category was 348 (Saharon Shelah).  Thus, though most authors post only a few papers in any given category, a few very prolific authors may distort our measurements of gender representation by paper-author.  Of course, these prolific people may also be very prominent within their subfields.  (Note that the most prolific author in our data set, Saharon Shelah, would not distort our measurement of gender representation, because his name was not automatically classified as either male or female.)  We graph the number of authors with different numbers of papers in Figure~\ref{F:papernumbers}.

In Table~\ref{Ta:femaleauthors}, we rank fields by the fraction of female authors out of all authors whose gender is known.  Again, we compute 95\% confidence intervals using a bootstrap method.

\begin{table}[!htbp]
\centering
\begin{tabular}{|c|c|c|}
\hline
\textbf{Subject} & \textbf{Fraction of Female Authors} & \textbf{Confidence Interval} \\
\hline

General Mathematics         &   0.088 & (0.062, 0.114) \\
Operator Algebras            &  0.111 & (0.095, 0.127) \\
Logic                        &  0.116 & (0.100, 0.133) \\
KTheory and Homology         &  0.116 & (0.092, 0.139) \\
Quantum Algebra              &  0.120 & (0.106, 0.134) \\
Category Theory              &  0.121 & (0.096, 0.145) \\
Complex Variables            &  0.127 & (0.112, 0.143) \\
Number Theory                &  0.132 & (0.122, 0.142) \\
Spectral Theory              &  0.133 & (0.113, 0.152) \\
Metric Geometry              &  0.135 & (0.118, 0.153) \\
Optimization and Control     &  0.135 & (0.124, 0.147) \\
Algebraic Topology           &  0.139 & (0.123, 0.156) \\
History and Overview         &  0.140 & (0.112, 0.167) \\
Differential Geometry        &  0.144 & (0.134, 0.154) \\
Symplectic Geometry          &  0.145 & (0.124, 0.168) \\
Numerical Analysis           &  0.149 & (0.137, 0.161) \\
Probability                  &  0.149 & (0.141, 0.158) \\
Group Theory                 &  0.150 & (0.137, 0.164) \\
Dynamical Systems            &  0.154 & (0.143, 0.165) \\
Functional Analysis          &  0.155 & (0.143, 0.167) \\
Algebraic Geometry           &  0.158 & (0.148, 0.168) \\
Classical Analysis and ODEs  &  0.164 & (0.150, 0.178) \\
Representation Theory        &  0.168 & (0.153, 0.182) \\
Geometric Topology           &  0.170 & (0.156, 0.184) \\
General Topology             &  0.174 & (0.143, 0.203) \\
Statistics Theory            &  0.176 & (0.163, 0.189) \\
Combinatorics                &  0.178 & (0.169, 0.187) \\
Analysis of PDEs             &  0.186 & (0.176, 0.196) \\
Rings and Algebras           &  0.190 & (0.173, 0.207) \\
Commutative Algebra          &  0.212 & (0.192, 0.232) \\
\hline

\end{tabular}
\caption{Fraction of Female Authors}\label{Ta:femaleauthors}
\end{table}

In Table~\ref{Ta:femaleauthorconfidence}, we use the confidence intervals for fraction of female authors to divide fields into categories, based on whether the confidence interval contains the average fraction of female authors, 0.155.  We see that this method of analysis places more subfields in the ``Average'' category than our analysis by paper-author; presumably, a few prolific people had an outsize impact on those subfields.

\begin{table}[!htbp]
\centering
\begin{tabular}{|c|p{10cm}|}
\hline
\textbf{Female Authors} & \textbf{\arxiv Categories} \\
\hline

Significantly Fewer & Category Theory, Complex Variables, Differential Geometry, General Mathematics, K-Theory and Homology, Logic, Metric Geometry, Number Theory, Operator Algebras, Optimization and Control, Probability, Quantum Algebra, Spectral Theory \\

\hline

Average & Algebraic Geometry, Algebraic Topology, Classical Analysis and ODEs, Dynamical Systems, Functional Analysis, General Topology, Group Theory, History and Overview, Numerical Analysis, Probability, Representation Theory, Symplectic Geometry \\

\hline

Significantly More & Analysis of PDEs, Combinatorics, Commutative Algebra, Geometric Topology, Rings and Algebras, Statistics Theory \\

\hline
\end{tabular}
\caption{Female Authors vs. Expectation}\label{Ta:femaleauthorconfidence}
\end{table}

\clearpage

\section{Questions and Hypotheses}

Our analysis shows that women are underrepresented as authors of mathematics papers on the \arxiv\!, even in comparison to the proportion of women who hold full-time positions in mathematics departments.  Moreover, our measurements demonstrate that there are significant differences in women's representation between different subfields of mathematics.  These measurements corroborate the anecdotal reports of differences between mathematical subfields.

A key feature of our analysis is our computation of confidence intervals using a bootstrap method.  Standard measures of statistical significance, such as the chi-squared test, are not always useful when applied to very large data sets.  In our case, for example, the chi-squared test tells us that every \arxiv category except for Dynamical Systems, General Topology, and Symplectic Geometry has a statistically significant different number of women paper-authors from the overall \arxiv average with a confidence level of $p<.02$.  For Algebraic Geometry, to pick an arbitrary category, we have $p<10^{-12}$.  However, algebraic geometry papers dominate our data set: we are not nearly as convinced that algebraic geometry is unusual as such a $p$-value might imply.  By grouping subfields where confidence intervals overlap, we obtain a more robust characterization of the similarities and differences between branches of mathematics.

In many cases, our measurements of representation in specific subfields are consistent with other measurements of women's contributions to these subfields.  Rates of \arxiv postings often correlate with PhD topics.  For example, the low rate of women PhDs in analysis matches our low measurements of women posting papers in Operator Algebras (despite the name, an analytic topic), and given the comparatively small number of women PhDs in physics, we might not be surprised to find few women posting in physics-influenced categories such as Quantum Algebra.  Conversely, comparatively high rates of women PhDs in combinatorics and statistics are matched by high rates of women's \arxiv postings in those areas.  Though the classification schemes are different, we also find broad agreement with the measurements of women's publications by subfield in \cite{zbmath}.

Using a finer classification scheme, we identify subfields with especially high or low representation of women whose properties are obscured in the AMS PhD data.  For example, Commutative Algebra and Combinatorics are two subfields with comparatively very high numbers of women publishing.  Sometimes a category that appears average in the PhD data actually combines subfields with widely differing participation of women: for instance, the Geometric Topology \arxiv category (where knot theory papers are posted) has among the highest rates of women while Differential Geometry has among the lowest, but both would be classified as Geometry/Topology in the AMS PhD report.

Understanding the representation of women in different subfields of mathematics helps us to evaluate claims that specific organizations have been particularly effective or ineffective at recruiting women participants.  For example, knowing the proportion of women in a particular research area allows us to judge whether a conference has made unusually successful efforts to recruit women speakers.  To take a representative case, Kristin Lauter wrote in the May-June 2016 President's Report for the Association for Women in Mathematics \cite{awmnewsletter}:

\begin{quotation}
Two long-running international biannual conferences in number theory
will take place this summer featuring many more women as plenary and invited
speakers than ever before: the 14th Meeting of the Canadian Number Theory
Association (CNTA XIV) has 3/7 female plenary speakers and 7/21 female invited
speakers; the Twelfth Algorithmic Number Theory Symposium (ANTS-XII) has 2/5
female invited speakers.
\end{quotation}

\noindent If we know only that women receive about a third of mathematics PhDs, these statistics look rather ordinary; when one realizes that number theory has below-average participation of women for a mathematics research field, a 30\% or 40\% rate of invited speakers seems more impressive.  (Note that a truly unbiased selection of conference speakers is likely to result in a \emph{higher} number of women than the population average; see \cite{martin} for analysis of some specific mathematics conferences, and \cite{calculator} for a simulation.)

Similarly, the authors of \cite{topazsen} studied the representation of women on the editorial boards of mathematical journals.  They found that the proportion of women on mathematics journal editorial boards is significantly lower than the proportion of women with faculty positions in mathematics at doctoral-granting institutions.  Since one would expect journals focused on subfields with large numbers of women to have more women on their editorial boards, understanding the participation of women in different subfields is important for evaluating a journal's record.  For example, Topaz and Sen write, ``Within our data set, the journals published by SIAM Publications have amongst the
highest representation of women.''  Because SIAM is a professional society for applied mathematicians, one might expect its editorial boards to reflect the greater participation of women in applied mathematics.  The observation is more surprising when one realizes that even SIAM journals for low-participation subfields, such as \emph{SIAM Journal on Optimization}, have an above-average fraction of women on their editorial boards.  (For some reactions to Topaz and Sen's investigation from SIAM editors, see \cite{BW}.)

Though useful, our measurements are limited by our use of given names to determine gender.  Some of the apparent underrepresentation of women on the \arxiv may be accounted for by the possibility that women are more likely than men to use only their initials on papers.  Because the use of initials varies by subfield, such a strategy may be more common in some \arxiv categories than others.  Differences in the gendering of given names across cultures, or in the names categorized in our reference data set, could introduce biases that over- or under-estimate the number of women from particular nations.  Although it is expensive to do so at scale, these limitations could be reduced by researching the gender of individual mathematicians.  We note that our classification also obscures the contributions of mathematicians whose gender identification is nonbinary.

Another limitation of this data set is the use of \texttt{arXiv} postings as a measure of a combination of representation and productivity.  As discussed in Section 2, this may underestimate the contributions of mathematicians who choose not to post preprints to the \texttt{arXiv}, whether because they make publications available in other ways, or because they devote more of their time and energy to teaching or service than to research.  We anticipate that differences in popularity of the \texttt{arXiv} among different fields of mathematics will affect men and women similarly, and thus have little impact on our comparisons of women's representation in different subfields.  We also note that our results are most reflective of underrepresentation of women mathematicians in positions involving research productivity; other data sets, such as abstracts or lists of attendees at conferences of the Mathematical Association of America or  American Mathematical Association of Two-Year Colleges, would be better suited to investigating gender distributions in primarily or exclusively teaching-oriented positions.

Our choice of classification scheme for subfields is of necessity somewhat arbitrary.  A single person's research may cross the lines between different \arxiv categories; on the other hand, somebody might identify with a subfield more specific than \arxiv categories can reflect.  Furthermore, the lines between subfields change over time, as research progresses and fashions change.  Our measure takes a snapshot of approximately twenty years, but is biased toward more recent time periods, since the rate of postings to the \arxiv has been increasing.  As we noted in \S~\ref{S:arxivintro}, we have used posted preprints, a form of publication, as a proxy for membership in subfields.  Thus, we are more likely to detect contributions from mathematicians who have the resources to publish and who participate in social structures that reward sharing preprints.

Differences in women's representation within mathematics warrant further investigation, as they may provide clues for increasing women's representation in the field, on the \arxiv and beyond.  What factors might affect the participation of women in specific mathematical subfields?

Many mathematicians favor the mentorship hypothesis: a few good mentors can have a strong positive effect on the participation of women in a subfield.  (Conversely, sexist or discriminatory actions by prominent people in a particular subfield might drive women away.)  For example, Judy Green and Jeanne LaDuke argue in \cite{pioneer} that a handful of supportive advisors made algebraic geometry the most popular dissertation field for women completing Ph.D.s in mathematics before 1940.

Other mathematicians doubt that a problem exists.  One possible response to observations about disproportionate representation within a field is to wonder whether there is some quality inherent to the overrepresented group which makes members of that group naturally more interested in or talented at that subject.  In our context, for example, one might advocate the hypothesis that men are simply more interested in operator algebras.  This proposal is belied by the varying proportions of women in STEM fields during different time periods and in different nations.  Hyde and Mertz \cite{HY09} found that the percentage of girls on International Mathematical Olympiad teams is significantly correlated with countries' Gender Gap Index, a measure of the difference in opportunities available to women and men. In our analysis, we found that the proportion of women paper-authors on the arXiv has increased over time (Figure~\ref{F:womenovertime}).  At every point in the past, we would have been mistaken to suppose that we had reached an upper bound on women's interest and ability in STEM; it seems egocentric for us as a society to suppose we have reached such a point now.  Further, the idea that men are inherently better suited than women to a particular subfield of math is not a testable hypothesis.  Other hypotheses, such as ``Implementing a blind application review process will increase the proportion of women who are hired,'' can be tested and affirmed or refuted \cite{AS07, BEH14}.  We suggest that practices currently in place in the mathematical subfields we identified with above-average representation of women, such as geometric topology, combinatorics, and commutative algebra, will provide fruitful ground for the generation of testable hypotheses.

Differences between nations may affect the representation of women in mathematics: different countries have different rates of women's participation in mathematics, and mathematical traditions focused on different subfields.  For example, the percentage of mathematicians in Italy who are female is higher than the percentage of mathematicians in the US who are female: 2005 data showed that 35\% of academic mathematicians in Italy were women \cite{hk}.  Meanwhile, Italy has a very strong tradition of research in algebraic geometry, dating back to the ``Italian school'' of the early twentieth century.  Thus, we might expect to find comparatively high numbers of Italian women in algebraic geometry.

Even within pure mathematics, applications may play a role.  We hypothesize that women are more likely to work on problems that have applications to other fields they find interesting.  Therefore, on average, we should find more women working on problems with applications to biology than on problems with applications to physics or computer science.  Subfields with plentiful or well-advertised applications to biology would then attract more women.

Accessibility of problems may also be important.  Many women in the United States first consider academic careers in mathematics after successful undergraduate research experiences.  If these women go on to specialize in similar fields, then subfields with many problems accessible to undergraduates might attract more women.

Finally, the authors of \cite{LCMF} showed that fields whose practitioners believe you must be brilliant to succeed have lower representation of women.  Furthermore, this predicts the low representation of women in mathematics as compared to women in statistics.  One might ask whether similar factors affect the participation of women in mathematical subfields: are some subfields of mathematics viewed as more dependent on unique insight?  In this case, an effort to make problems in such subfields more accessible might draw more women into mathematics.

%Questions about women's participation in mathematics pose many intriguing research and policy questions.

\end{document}